\documentclass{amsart}
\usepackage{amsthm}                 %Definitionen, Theoreme etc.
\usepackage{amsmath}                %siehe userguide, z.B. boxed{}, text-Befehl etc.
\usepackage{amsfonts,amssymb}       %ein paar zusätzliche Symbole
\usepackage{graphicx}               %Bilder
\usepackage{float}                  %Bilder/Positionierung
\usepackage{rotating}               %Bilder/Positionierung
\restylefloat{figure}               %Bilder/Positionierung
\usepackage{mathrsfs}               %mathscr-Buchstaben
\usepackage{fancyhdr}               %Kopfzeile
\usepackage[breaklinks]{hyperref}   %Referenzen anklicken, Zeilenumbruch im Inhaltsverzeichnis
\usepackage[numbers]{natbib}
\theoremstyle{plain}
\newtheorem{lem}{Lemma}

\newtheorem{thm}[lem]{Theorem}
\newtheorem{cor}[lem]{Corollary}

\theoremstyle{definition}

\newtheorem{rem}[lem]{Remark}
\newcommand{\R}{\mathbb R}

\newcommand{\Z}{\mathbb Z}

\newcommand{\Diff}{\mbox{\rm Diff}}

\renewcommand{\L}{\mathcal L}
\newcommand{\A}{\mathbb A}
\newcommand{\id}{\text{\rm id}}
\newcommand{\dx}{\,\text{\rm d}x}
\renewcommand{\d}{\,\text{\rm d}}
\newcommand{\dw}{\text{\rm d}}

\newcommand{\g}{\mathfrak{g}}

\newcommand{\X}{\mathfrak{X}}

\renewcommand{\S}{\mathbb S}

\renewcommand{\phi}{\varphi}

\newcommand{\norm}[1]{\left|\!\left|#1\right|\!\right|}

\newcommand{\eps}{\varepsilon}

\newcommand{\ska}[2]{\left\langle #1,#2\right\rangle}

\newcommand{\set}[2]{\left\{#1;\;#2\right\}}

\newcommand{\bea}{\begin{eqnarray}}
\newcommand{\eea}{\end{eqnarray}}
\newcommand{\beq}{\begin{equation}}
\newcommand{\eeq}{\end{equation}}
\renewcommand{\phi}{\varphi}

\renewcommand{\autoref}[1]{\text{Eq.}~\eqref{#1}}
\begin{document}
\title{On a two-component $\pi$-Camassa--Holm system}
\author{Martin Kohlmann}
\address{Peter L. Reichertz Institute for Medical Informatics, University of Braunschweig, D-38106 Braunschweig, Germany}
\email{martin.kohlmann@plri.de}
\begin{abstract} A novel $\pi$-Camassa--Holm system is studied as a geodesic flow on a semidirect product obtained from the diffeomorphism group of the circle. We present the corresponding details of the geometric formalism for metric Euler equations on infinite-dimensional Lie groups and compare our results to what has already been obtained for the usual two-component Camassa--Holm equation. Our approach results in well-posedness theorems and explicit computations of the sectional curvature.
\end{abstract}
\maketitle
\tableofcontents
\section{Introduction}\label{sec_intro}
In this paper, we are concerned with the following variation of the two-component Camassa-Holm system:
\bea
\left\{
\begin{array}{rcl}
m_t&=&-m_xu-2u_xm-\pi(\rho)\rho_x,\\
\pi(\rho)_t&=&-(\pi(\rho)u)_x
\end{array}
\right.
\label{2piCH}\eea
where
$$m=u-u_{xx}\quad\text{and}\quad\pi(\rho)=\rho-\mu(\rho)=\rho-\int_0^1\rho\dx;$$
here $u(t,x)$ and $\rho(t,x)$ depend on a time variable $t\geq 0$ and a space variable $x\in\S=\R/\Z$. We will call the system \eqref{2piCH} the $\pi$-2CH equation. Observe that \autoref{2piCH} reduces for $\mu(\rho)=0$ to the two-dimensional Camassa-Holm (2CH) equation studied in, e.g., \cite{CLZ06,CI08,F06,HT09}, and for $\pi(\rho)=0$ to the one-component Camassa-Holm (CH) equation \cite{CH93}.

The CH equation first appeared in 1993 and was introduced by virtue of a bi-Hamiltonian approach in \cite{CH93}. In the subsequent years, the properties of its solutions, which exhibit typically nonlinear phenomena as wave breaking or peakons, have been examined in detail, \cite{CE00,CM00,CS02,M02,M04}. In 1998, Misio{\l}ek \cite{M98} showed that the CH equation re-expresses geodesic motion on the Bott-Virasoro group. A little later, the CH equation with periodic boundary conditions had been recast as geodesic flow on the diffeomorphism group of the circle $\S$, \cite{K99,CK02,CK03}. A condensed account of the geometric picture for the CH equation can be found in \cite{K04,L07}. Similarly, the 2CH equation has been discussed under geometric aspects and the most important results can be found in \cite{EKL11,HT09}; here it is shown that 2CH corresponds to a geodesic flow on a semidirect product Lie group which is obtained from the diffeomorphism group of the circle. Some qualitative properties of the solutions of 2CH and integrability issues are the subject of \cite{CI08,ELY07,GZ10,LL09}.

We will consider the system \eqref{2piCH} in the spaces $H^s\times H^{s-1}/\R$ for $s>5/2$, where $H^s=H^s(\S)$ denotes the $L_2$-Sobolev space of regularity $s$ on the circle. There are three main reasons why a study of the system \eqref{2piCH} is of interest:
\begin{itemize}
\item First, since for the 2CH system
$$\frac{\dw}{\dw t}\int_\S\rho\dx=-\int_\S(\rho u)_x\dx=0,$$
it follows that if $\rho$ has mean zero at time $t=0$, the mean of the solution $\rho$ to 2CH will have mean zero at all later times.
\item Second, we may decompose the space $H^{s-1}=\hat H^{s-1}\oplus\R$, where $\hat H^{s-1}$ is the subspace of $H^{s-1}$ containing all zero mean functions, with the associated projections
    $$(\pi,\mu)\colon H^{s-1}\to\hat H^{s-1}\times\R,\quad\rho\mapsto(\pi(\rho),\mu(\rho)),$$
    so that $\rho=\pi(\rho)+\mu(\rho)$ for any $\rho\in H^{s-1}$. The aim of the present paper is to comment on the behavior of the 2CH system on the non-trivial component $H^s\times\hat H^{s-1}$. The interesting aspect is that 2CH has a meaningful geometric interpretation on the entire space $H^s\times H^{s-1}$ as well as on the component $H^s\times\hat H^{s-1}$ which will be provided by our results.
\item Finally, a related $\pi$-2HS system has been proposed in \cite{L11},
\bea
\left\{
\begin{array}{rcl}
u_{txx}&=&-2u_xu_{xx}-uu_{xxx}+\pi(\rho)\rho_x,\\
\pi(\rho)_t&=&-(\pi(\rho)u)_x
\end{array}
\right.
\nonumber
\eea
where the author establishes its equivalence to a geodesic equation on a K\"ahlerian manifold which is isometric to a subset of a complex projective space, providing an example of an infinite-dimensional Hopf fibration.
\end{itemize}
The outline of the paper and its main results are as follows: In Section~\ref{sec_prelim} we recall some preliminaries, concerning diffeomorphism groups, semidirect products and Arnold's geometric formalism \cite{A66} which will be applied to \autoref{2piCH} in the following. In Section~\ref{sec_geomsol} we show that \autoref{2piCH} can be recast as a geodesic equation on a suitable semidirect product group. From the geometric formulation, we obtain well-posedness of the system \eqref{2piCH} in $H^s\times H^{s-1}/\R$ and in $C^{\infty}\times C^{\infty}/\R$.
Finally, in Section~\ref{sec_curv} we perform a lengthy computation of the sectional curvature $S$ associated with the group of $\pi$-2CH and compare it to corresponding results for the usual two-component Camassa-Holm system.
\section{Preliminaries}\label{sec_prelim}
As brief as possible, since the material appears in many publications dealing with geometric aspects of the Camassa-Holm equations, we provide some background information which is necessary for this paper.

\subsection{The diffeomorphism group of the circle}
Let $\S=\R/\Z$ and $s\geq 0$. We denote by $H^s=H^s(\S)$ the $L_2$-Sobolev space of order $s$ on the circle. Let $H^s\Diff(\S)$ denote the set of orientation-preserving diffeomorphisms $\S\to\S$ in $H^s$. It is well-known that $H^s\Diff(\S)$ is a topological group (with respect to composition) and a smooth Hilbert manifold for any $s>3/2$, cf.~\cite{EM70}; an atlas is given by the charts $(U_i,\Phi_i)$, $i=1,2$, where
\bea U_1&=&\set{u\in H^s}{u_x>-1,-\tfrac{1}{2}<u(0)<\tfrac{1}{2}},\nonumber\\
U_2&=&\set{u\in H^s}{u_x>-1,0<u(0)<1}\nonumber\eea
and
$$\Phi_i\colon U_i\to H^s\Diff(\S),\quad \Phi_i(u)=\id+u,$$
for $i=1,2$, cf.~\cite{GR05}. The tangent space of $H^s\Diff(\S)$ at the identity can be identified with the $H^s$ vector fields on the circle and hence with $H^s$. Furthermore $H^s\Diff(\S)$ is parallelizable, i.e., one has the trivialization
$$TH^s\Diff(\S)\simeq H^s\Diff(\S)\times H^s,$$
and the derivative of the right translation map $R_\phi\colon\psi\mapsto\psi\circ\phi$ on $H^s\Diff(\S)$ is an automorphism of $H^s$. For $s\to\infty$, the groups $H^s\Diff(\S)$ approximate the Lie group $C^{\infty}\Diff(\S)$ of smooth and orientation-preserving diffeomorphisms $\S\to\S$, a $C^{\infty}$-Fr\'echet manifold on which inversion and composition are smooth maps. We can describe its manifold structure as above by replacing $H^s$ with $C^{\infty}(\S)=C^{\infty}$. Note that $T_\id C^{\infty}\Diff(\S)\simeq C^{\infty}$ is a Lie algebra with the bracket
$$[u,v]=v_xu-u_xv.$$
\subsection{Semidirect products}
Let $G$ be a Lie group and $V$ be a vector space. If $G$ acts on
the right on $V$, one defines
$$(g_1,v_1)(g_2,v_2)=(g_1g_2,v_2+v_1g_2)$$
and with this product, $G\times V$ becomes a Lie group (the \emph{semidirect product of $G$ and $V$}) which is denoted as $G\circledS V$. It is easy to see
that $(e,0)$ is the neutral element, where $e$ denotes the neutral
element of $G$, and that $(g,v)$ has the inverse
$(g^{-1},-vg^{-1})$. The Lie bracket on the Lie algebra of $G\circledS V$ is given by
$$[(\xi_1,v_1),(\xi_2,v_2)]=([\xi_1,\xi_2],v_2\xi_1-v_1\xi_2),$$
where $v\xi$ denotes the induced action of the Lie algebra $\g$ on $V$ and $[\cdot,\cdot]$ is the Lie bracket on $\g$.

We denote by $H^{s-1}/\R$ the space $H^{s-1}$ with two functions being identified if they differ by a constant, and write $[\rho]$ for the elements of $H^{s-1}/\R$. We let $G^s$ be the semidirect product $H^s\Diff(\S)\circledS H^{s-1}(\S)/\R$ and $G^{\infty}=C^{\infty}\Diff(\S)\circledS C^{\infty}(\S)/\R$. The
group product in these groups is given by
$$(\phi_1,[f_1])(\phi_2,[f_2]):=(\phi_1\circ\phi_2,[f_2+f_1\phi_2])$$
where $f\phi:= f\circ\phi$ is a right action of $H^s\Diff(\S)$ on the scalar functions on $\S$. The neutral element is $(\id,[0])$ and $(\phi,[f])$ has the inverse $(\phi^{-1},-[f\circ\phi^{-1}])$. Clearly, $G^s$ is a smooth Hilbert manifold and a topological group and $G^{\infty}$ is a smooth Fr\'echet manifold and a Lie group. We have the trivializations
$$TG^s\simeq G^s\times H^s\times H^{s-1}/\R,\quad TG^\infty\simeq G^\infty\times C^{\infty}\times C^{\infty}/\R$$
and the Lie bracket on $C^{\infty}\times C^{\infty}/\R$ is given by
$$[(u_1,[u_2]),(v_1,[v_2])]=([u_1,v_1],[v_{2x}u_1-u_{2x}v_1]).$$
The derivative of the right shift operator $R_{(\phi,[f])}\colon(\psi,[g])\mapsto(\psi,[g])(\phi,[f])$ on these groups is
\beq\label{rightshift}DR_{(\phi,[f])}(v,[h])=(v,[h])\circ\phi\eeq
and it is an automorphism of $H^s\times H^{s-1}/\R$ or $C^{\infty}\times C^{\infty}/\R$ respectively.
\subsection{Euler equations on infinite-dimensional Lie groups} Consider a rigid body in $\R^3$ with three rotational degrees of freedom. The classical Euler equation for the motion of the body can also be interpreted in a geometric framework, i.e., as geodesic equation on the finite-dimensional Lie group $SO(3)$. The geodesics are length-minimizing with respect to the left-invariant metric on $SO(3)$ that is defined by the inertia matrix of the body. The geodesic equation can be written down in terms of a set of Christoffel symbols, as a second order equation
$$\dot u^k+\Gamma^k_{ij}u^iu^j=0$$
for the Eulerian velocity $u$ of the rigid body. The Eulerian velocity is obtained from the Lagrangian picture by applying the derivative of the left shift operator on $SO(3)$ to the velocity on $SO(3)$. The Christoffel symbols in turn define a covariant derivative which is compatible with the left-invariant metric.

According to Arnold's \cite{A66} fundamental observation and the work of Ebin and Marsden \cite{EM70}, the above formalism works analogously for the Camassa-Holm equation and its supersymmetric extension which have the configuration manifolds $C^{\infty}\Diff(\S)$ and $C^{\infty}\Diff(\S)\circledS C^{\infty}$ respectively. However, the passage from the finite-dimensional group $SO(3)$ to the infinite-dimensional diffeomorphism groups has to be carried out with care: The kinetic energy metric on the diffeomorphisms has to be right-invariant instead of left-invariant in order to obtain the correct equations of motion. It is induced by an inertia operator which is a topological isomorphism of the Lie algebra. The Christoffel symbols turn into a Christoffel operator which is a right-invariant bilinear map on the tangent bundle. The Eulerian velocity is now obtained from a right shift of the Lagrangian velocity to the Lie algebra. Furthermore, instead of using the Riemannian geometry of finite-dimensional manifolds we now have to apply the Riemannian geometry for general Banach manifolds which is more subtle. Instead of presenting more details for CH and 2CH, we refer the reader to \cite{EKL11,HMR98,HT09,K04} for an extensive presentation and use the $\pi$-2CH equation to exemplify the approach once again in our next section.
\section{The geometry and solutions of the $\pi$-2CH system}\label{sec_geomsol}
In this section, we will straightforwardly introduce the geometric picture for \autoref{2piCH} from which we immediately obtain some well-posedness results. For technical purposes, the configuration manifold for $\pi$-2CH will be $G^s$ at first.
\subsection{The inertia operator}
Fix $s>5/2$. We introduce $A=1-\partial_x^2\colon H^{s}\to H^{s-2}$, $B=\pi\colon H^{s-1}/\R\to\hat H^{s-1}$ and $\A=\text{diag}(A,B)$ with domain $D(\A)=H^{s}\times H^{s-1}/\R$ and range $H^{s-2}\times\hat H^{s-1}$. It is easy to see that $\A$ is well-defined and a topological isomorphism.
\subsection{The right-invariant metric} We introduce the bilinear and symmetric map
$$\ska{\cdot}{\cdot}\colon (H^{s}\times H^{s-1}/\R)^2\to\R,\quad((u,[\rho]),(v,[\tau]))\mapsto\int_\S(u,[\rho])\A(v,[\tau])\dx$$
which is explicitly given as
$$\ska{(u,[\rho])}{(v,[\tau])}=\int_\S(uv+u_xv_x+\rho\tau)\dx-\mu(\rho)\mu(\tau).$$
To see that $\ska{\cdot}{\cdot}$ is positive definite, we observe that
$$\ska{(u,[\rho])}{(u,[\rho])}=\norm{u}_{H^1}^2+\mu(\rho^2)-\mu(\rho)^2$$
and that $\mu(\rho)^2\leq\mu(\rho^2)$ by the Cauchy-Schwarz inequality, with equality if and only if $1$ and $\rho$ are linearly dependent, i.e., iff $[\rho]=[0]$. Since $H^s\times H^{s-1}/\R\simeq T_\id G^s$ and in view of the identity \eqref{rightshift} we now extend $\ska{\cdot}{\cdot}$ to a (weak) right-invariant metric on $G^s$ by setting
$$\ska{(U_1,[U_2])}{(V_1,[V_2])}_{(\phi,[f])}=\ska{(U_1,[U_2])\circ\phi^{-1}}{(V_1,[V_2])\circ\phi^{-1}},$$
for all $(U_1,[U_2])(V_1,[V_2])\in T_{(\phi,[f])}G^s\simeq H^s\times H^{s-1}/\R$. That the map $(\phi,[f])\mapsto\ska{\cdot}{\cdot}_{(\phi,[f])}\in\L^2_{\text{sym}}(H^s\times H^{s-1}/\R;\R)$\footnote{For Banach spaces $X$ and $Y$, $\L^2_{\text{sym}}(X;Y)$ is the space of symmetric bilinear maps $X\to Y$.} is indeed smooth follows from the representation
\bea\ska{(U_1,[U_2])}{(V_1,[V_2])}_{(\phi,[f])}&=&\int_S\left((U_1V_1+U_2V_2)\phi_x+\frac{U_{1x}V_{1x}}{\phi_x}\right)\dw x\nonumber\\
&&\quad-\int_\S U_2\phi_x\dx\int_\S V_2\phi_x\dx,\nonumber\eea
and that the natural topology on any $T_{(\phi,[f])}G^s$ is stronger than the topology induced by $\ska{\cdot}{\cdot}_{(\phi,[f])}$ follows as in Par.~9 of \cite{EM70}.
\subsection{The Christoffel operator}
To obtain the Christoffel map for \autoref{2piCH}, we rewrite \autoref{2piCH} as
\beq
\left(
  \begin{array}{c}
    u_t+uu_x \\
    \left[\rho\right]_t+\left[u\rho_x\right] \\
  \end{array}
\right)
=
\left(
  \begin{array}{c}
    -\frac{1}{2}A^{-1}\partial_x(2u^2+u_x^2+\pi(\rho)^2) \\
    -[u_x\pi(\rho)] \\
  \end{array}
\right);
\label{2piCHrew}
\eeq
in this form, the system is also suitable for the formulation of weak solutions. That \eqref{2piCHrew} is indeed equivalent to the $\pi$-2CH equation in its initial form follows by applying the operator $\A$ to \autoref{2piCHrew}: The first component reproduces the equation with $m_t$ in \eqref{2piCH} and the second component gives
$$\pi(\rho)_t+\pi(u\rho_x)=-\pi(u_x\pi(\rho))$$
which can equivalently be written as
$$\pi(\rho)_t=-\pi((u\pi(\rho))_x)=-(\pi(\rho)u)_x.$$
We now define
$$
\Gamma(u,v)=-\frac{1}{2}
\left(
  \begin{array}{c}
    A^{-1}\partial_x(2u_1v_1+u_{1x}v_{1x}+\pi(u_2)\pi(v_2)) \\
    \left[u_{1x}\pi(v_2)+v_{1x}\pi(u_2)\right] \\
  \end{array}
\right),
$$
for all $(u,v)\in H^s\times H^{s-1}/\R\simeq T_{(\id,[0])}G^s$, and extend $\Gamma$ to a right-invariant map $\Gamma_{(\phi,[f])}$ on $TG^s$ by setting
$$\Gamma_{(\phi,[f])}(U,V)=\Gamma(U\circ\phi^{-1},V\circ\phi^{-1})\circ\phi,\quad U,V\in T_{(\phi,[f])}G^s\simeq H^s\times H^{s-1}/\R.$$
Then $\Gamma$ defines a smooth spray $(\phi,[f])\mapsto\Gamma_{(\phi,[f])}$, $G^s\to\mathcal L_{\text{sym}}^2(H^s\times H^{s-1}/\R;H^s\times H^{s-1}/\R)$, see \cite{EKL11} for a proof in a similar situation.
\subsection{The torsion-free affine connection}
Let $\X$ denote the space of smooth vector fields on $G^s$ and define, for any $X,Y\in\X$ the map
$$\nabla_XY=DY\cdot X-\Gamma(X,Y),$$
where $DY(\phi,[f])\cdot X(\phi,[f])=\left.\frac{\dw}{\dw\eps}\right|_{\eps=0}Y((\phi,[f])+\eps X(\phi,[f]))$. Then $\nabla$ defines a torsion-free affine connection on $G^s$, i.e.,
\begin{itemize}
\item[(i)] $\nabla_{fX+gY}Z=f\nabla_XZ+g\nabla_YZ$,
\item[(ii)] $\nabla_X(Y+Z)=\nabla_XY+\nabla_XZ$,
\item[(iii)] $\nabla_X(fY)=f\nabla_XY+X(f)Y$,
\item[(iv)] $\nabla_XY-\nabla_YX=[X,Y]$,
\end{itemize}
for all $X,Y,Z\in\mathfrak X$ and all $f,g\in C^{\infty}(G^s;\R)$; this can be proved using that the Lie bracket is given locally by
$$[X,Y]=DY\cdot X-DX\cdot Y.$$
\subsection{The geodesic flow and well-posedness}
A geodesic on $G^s$ is a solution $t\mapsto (\phi,[f])(t)$ of the equation $\nabla_{(\phi_t,[f]_t)}(\phi_t,[f]_t)=0$, i.e.,
$$(\phi_{tt},[f]_{tt})=\Gamma_{(\phi,[f])}((\phi_t,[f]_t),(\phi_t,[f]_t)).$$
On the other hand, since we have a metric on $G^s$, a geodesic can also be understood as a minimizer of the functional
$$L(\gamma)=\int_J\ska{\gamma_t(t)}{\gamma_t(t)}_{\gamma(t)}^{1/2}\d t,\quad\gamma\in C^1(J,G^s).$$
We now show that the metric and the affine connection produce indeed the same geodesic flow on $G^s$, i.e., the connection $\nabla$ preserves the metric $\ska{\cdot}{\cdot}$.
\lem The metric $\ska{\cdot}{\cdot}$ and the connection $\nabla$ are compatible in the usual sense
$$X\ska{Y}{Z}=\ska{\nabla_XY}{Z}+\ska{Y}{\nabla_XZ},\quad\forall X,Y,Z\in\mathfrak X.$$
\endlem\rm
\proof By the same arguments as in the proof of Proposition~3.1 in \cite{EKL11} it is enough to check that
\bea 0&=&\ska{v_{2x}u_1}{Bw_2}+\ska{w_{2x}u_1}{Bv_2}+\ska{\tfrac{1}{2}(\pi(v_2)\pi(u_2))_x}{w_1}\nonumber\\
&&+\ska{\tfrac{1}{2}(\pi(w_2)\pi(u_2))_x}{v_1}+\ska{\tfrac{1}{2}(v_{1x}\pi(u_2)+u_{1x}\pi(v_2))}{Bw_2}
\nonumber\\
&&+\ska{\tfrac{1}{2}(w_{1x}\pi(u_2)+u_{1x}\pi(w_2))}{Bv_2},\nonumber\eea
where $\ska{\cdot}{\cdot}$ denotes the $L_2$-pairing. Using the definition of $\mu$ and $\pi$ and performing integration by parts, we see that this is indeed true.
\endproof
We thus obtain a well-defined geodesic flow for the $\pi$-2CH equation. Since the geodesic spray for $\pi$-2CH is smooth, we also obtain the following well-posedness result.
\thm\label{thm_lwpgeo} Let $s>5/2$. There is an open neighborhood $U$ of zero in $H^s\times H^{s-1}/\R$ such that for any $(u_0,[\rho_0])\in U$ there is $T>0$ and a unique solution $(\phi,[f])\in C^{\infty}([0,T),G^s)$ of the initial value problem
\bea
\left\{
\begin{array}{rcl}
  (\phi_{tt},[f]_{tt}) & = & \Gamma_{(\phi,[f])}((\phi_t,[f]_t),(\phi_t,[f]_t)), \\
  (\phi_t,[f]_t)(0)    & = & (u_0,[\rho_0])  \\
  (\phi,[f])(0)        & = & (\id,[0])
\end{array}
\right.
\label{IVPgeo}
\eea
for the geodesic flow corresponding to the $\pi$-2CH system on $G^s$, with smooth dependence on $(u_0,[\rho_0])$, i.e., the local flow
$$\Phi\colon(t,u_0,[\rho_0])\mapsto \Phi(t,u_0,[\rho_0])=(\phi,[f]),\quad[0,T)\times U\to G^s$$
is smooth.
\endthm\rm
\cor Let $s>5/2$. There is an open neighborhood $U$ containing $(0,[0])$ in $H^s\times H^{s-1}/\R$ such that for any $(u_0,[\rho_0])\in U$ there is $T>0$ and a unique solution $(u,[\rho])$ to the initial value problem for the $\pi$-2CH system \eqref{2piCH} with
$$(u,[\rho])\in C([0,T);H^s\times H^{s-1}/\R)\cap C^1([0,T);H^{s-1}\times H^{s-2}/\R),$$
$(u,[\rho])(0)=(u_0,[\rho_0])$ and with continuous dependence on $(u_0,[\rho_0])$, i.e., the mapping
$$(u_0,[\rho_0])\mapsto (u,[\rho]),\quad U\to C([0,T);H^s\times H^{s-1}/\R)\cap C^1([0,T);H^{s-1}\times H^{s-2}/\R)$$
is continuous.
\endcor\rm
\proof Let $(\phi,[f])$ denote the solution of \autoref{IVPgeo} obtained in the proof of Theorem~\ref{thm_lwpgeo}. We now set
\beq(u,[\rho])=(\phi_t,[f]_t)\circ\phi^{-1}\label{sol}\eeq
and conclude, in view of the group properties of $G^s$, that $(u,[\rho])$ is a solution with the desired regularity to the Cauchy problem for the $\pi$-2CH equation.
\endproof
Indeed, well-posedness also holds in the smooth category. To prove this, we make use of the following lemma.
\lem\label{lem_conservation} The $\pi$-2CH equation enjoys the conservation laws
$$
(m\circ\phi)\phi_x^2+(\pi(\rho)\circ\phi) f_x\phi_x=m_0,\quad(\pi(\rho)\circ\phi)\phi_x=\pi(\rho_0).
$$
\endlem\rm
\proof Using the relation \eqref{sol} and \autoref{2piCH} we have that
$$\frac{\dw}{\dw t}[(\pi(\rho)\circ\phi)\phi_x]=[(\pi(\rho)_t+\pi(\rho)_xu+u_x\pi(\rho))\circ\phi]\phi_x=0$$
and second
$$\frac{\dw}{\dw t}[(m\circ\phi)\phi_x^2+(\pi(\rho)\circ\phi) f_x\phi_x]=[(m_t+m_xu+2u_xm+\pi(\rho)\rho_x)\circ\phi]\phi_x^2=0.$$
This achieves the proof.
\endproof
\thm There is an open neighborhood $U$ of zero in $H^3\times H^{2}/\R$ such that for any $(u_0,[\rho_0])\in U\cap (C^{\infty}\times C^{\infty}/\R)$ there is $T>0$ and a unique solution $(u,[\rho])$ to the initial value problem for the $\pi$-2CH system \eqref{2piCH} with
$$(u,[\rho])\in C^\infty([0,T);C^{\infty}\times C^{\infty}/\R),$$
$(u,[\rho])(0)=(u_0,[\rho_0])$ and with smooth dependence on $(u_0,[\rho_0])$, i.e., the map
$$(u_0,[\rho_0])\mapsto (u,[\rho]),\quad U\cap (C^{\infty}\times C^{\infty}/\R)\to C^\infty([0,T);C^{\infty}\times C^{\infty}/\R)$$
is smooth.
\endthm\rm
\proof We show that the geodesic flow obtained in Theorem~\ref{thm_lwpgeo} preserves spatial regularity when we increase the regularity of the initial data. Let $\Phi\colon[0,T_3)\times U_3\to G^3$, $\Phi(t,u_0,[\rho_0])=(\phi,[f])$, be the smooth local flow for the $\pi$-2CH equation on $G^3$, obtained in Theorem~\ref{thm_lwpgeo}. Pick $(u_0,[\rho_0])\in U_3\cap(C^{\infty}\times C^{\infty}/\R)$ and let $[0,T_s)$ denote the interval of existence for the associated solution $(\phi,[f])(t)$ in $G^s$, for $s\geq 3$. By uniqueness, we know that $T_s\leq T_3$ and it follows from Theorem 12.1 in \cite{EM70} and Lemma~\ref{lem_conservation} that $T_s<T_3$ is not possible. Hence $(\phi,[f])\in C^{\infty}([0,T_3);G^s)$ for any $s\geq 3$, and \autoref{sol}, the group properties of $G^{\infty}=\cap_{s>5/2} G_s$ and an application of Lemma 3.10 in \cite{EKK10} complete the proof.
\endproof
\section{The sectional curvature associated with the $\pi$-2CH system}\label{sec_curv}
For $X,Y,Z\in\X$ let
$$R(X,Y)Z=\nabla_X\nabla_YZ-\nabla_Y\nabla_XZ-\nabla_{[X,Y]}Z$$
and
$$S(X,Y)=\ska{R(X,Y)Y}{X}$$
denote the curvature tensor and the unnormalized sectional curvature associated with the $\pi$-2CH system. In this section, we compute an explicit formula for $S$ and compare our result to what has been obtained for the two-dimensional Camassa-Holm system in \cite{EKL11}. Note also that the physical motivation for studying curvatures is that the positivity (or negativity) of the sectional curvature of the configuration manifold is closely related to the stability (or instability) of the geodesic flow under small perturbations of the initial value \cite{A89}.
\thm\label{thm_curv} The sectional curvature $S$ at the identity element of $G^{\infty}$ for the $\pi$-2CH equation is given by
\bea
S(u,v)\!\!\!&=&\!\!\!\ska{\Gamma(u,v)}{\Gamma(u,v)}-\ska{\Gamma(u,u)}{\Gamma(v,v)}\nonumber\\
&&+\mu(u_{1x}v_2)^2+\mu(u_{2x}v_1)^2+\mu(u_1u_{2x})\mu(v_{1x}v_2)+\mu(u_2v_{1x})\mu(u_1v_{2x}).
\nonumber\eea
\endthm\rm
\proof We make use of the local formula for the tensor field $R_p$
\begin{align*}
R_p(U, V)W = & D_1\Gamma_{p}(W, U)V - D_1\Gamma_{p}(W, V)U
         \\
 &  + \Gamma_{p}(\Gamma_{p}(W, V), U) - \Gamma_{p}(\Gamma_{p}(W, U), V)
\end{align*}
where $D_1$ denotes differentiation with respect to $p$, i.e.,
$$D_1\Gamma_{p}(W, U)V = \frac{\dw}{\dw\eps}\bigg|_{\eps = 0} \Gamma_{p + \eps V}(W, U).$$
By right-invariance, it suffices to study the curvature at the identity $p=(\id,0)$.
Let $u = (u_1, u_2)$, $v = (v_1, v_2)$ and $w = (w_1, w_2)$ be three vectors in
$H^s\times H^{s-1}/\R$; to keep the notation as simple as possible, we omit the $[\cdot]$ for the second components in the following. First, we observe that
$$D_1\Gamma(w,u)v=-\Gamma(w_xv_1,u)-\Gamma(u_xv_1,w)+\Gamma(w,u)_xv_1.$$
Thus,
\begin{align}\nonumber
S(u,v) =&\; \ska{\Gamma(\Gamma(v,v),u)}{u}-\ska{\Gamma(\Gamma(v,u),v)}{u}
    \\ \label{Suv}
& + \langle \Gamma(v,u)_xv_1, u \rangle - \langle \Gamma(v,v)_xu_1, u \rangle
    \\ \nonumber
&+\ska{-\Gamma(v_xv_1,u)-\Gamma(v,u_xv_1)+2\Gamma(v_xu_1,v)}{u}.
\end{align}
To simplify the first four terms on the right-hand side of \autoref{Suv}, we define a bilinear operator $\mathcal B = (\mathcal B_1, \mathcal B_2)$ on $H^s\times H^{s-1}/\R$ by
\begin{align}
\begin{pmatrix}
  \mathcal B_1(u,v) \\
  \mathcal B_2(u,v)
\end{pmatrix}
=&\; \begin{pmatrix}
  -A^{-1}(2v_{1x}Au_1+v_1Au_{1x}+v_{2x}Bu_2) \\
  -B^{-1}(v_1Bu_2)_x
\end{pmatrix}
\nonumber
\\
=&\; \begin{pmatrix}
  -A^{-1}(2v_{1x}Au_1+v_1Au_{1x}+v_{2x}\pi(u_2)) \\
  -[(\pi(u_2)v_1)_x]
\end{pmatrix}.
\nonumber
\end{align}
Then $\mathcal B$ satisfies $\ska{\mathcal B(u,v)}{w}=\ska{u}{[v,w]}$ and
$$
\Gamma(u,v)=\frac{1}{2}\left[\begin{pmatrix}
  (u_1v_1)_x \\
  [u_{2x}v_1+v_{2x}u_1]
\end{pmatrix}
+B(u,v)+B(v,u)\right],
$$
and it follows by the line of arguments in Proposition~5.1 of \cite{EKL11} that we can simplify \eqref{Suv} to
\begin{align} \nonumber
S(u,v)=&\;\ska{\Gamma(u,v)}{\Gamma(u,v)}-\ska{\Gamma(u,u)}{\Gamma(v,v)}
    \\ \nonumber
&-\ska{\begin{pmatrix}
  u_{1x}v_1 \\
  u_{2x}v_1
\end{pmatrix}}{\Gamma(u,v)}+\ska{
\begin{pmatrix}
  u_{1x}u_1 \\
  u_{2x}u_1
\end{pmatrix}}{\Gamma(v,v)}
    \\ \nonumber
&+\ska{-\Gamma(v_xv_1,u)-\Gamma(v,u_xv_1)+2\Gamma(v_xu_1,v)}{u}.
\end{align}
It is a long and strenuous computation which shows that
\bea
&&
-\ska{\begin{pmatrix}
  u_{1x}v_1 \\
  u_{2x}v_1
\end{pmatrix}}{\Gamma(u,v)}+
\ska{
\begin{pmatrix}
  u_{1x}u_1 \\
  u_{2x}u_1
\end{pmatrix}}{\Gamma(v,v)} \nonumber\\
&&+\ska{-\Gamma(v_xv_1,u)-\Gamma(v,u_xv_1)+2\Gamma(v_xu_1,v)}{u}\nonumber\\
&=&\mu(u_{1x}v_2)^2+\mu(u_{2x}v_1)^2+\mu(u_1u_{2x})\mu(v_{1x}v_2)+\mu(u_2v_{1x})\mu(u_1v_{2x});
\nonumber
\eea
the details are left to the reader.
\endproof
\rem Note that a further consequence of Theorem~\ref{thm_curv} is that it provides another example that the sectional curvature is generally \emph{not} given by the nice formula
\beq\label{formulaS}S(u,v)=\ska{\Gamma(u,v)}{\Gamma(u,v)}-\ska{\Gamma(u,u)}{\Gamma(v,v)}\eeq
which turned out to be valid for the one-component Camassa-Holm equation and its supersymmetric extension (where $\ska{\cdot}{\cdot}$ and $\Gamma$ denote the corresponding metric and the Christoffel operator respectively). We have found a further example pointing out that an identity of type \eqref{formulaS} is far from being standard from the general theory of Riemannian connections on Banach manifolds. Other counterexamples are presented in \cite{KLM08,L11}.
\endrem\rm
\end{document}